\documentclass{amsart}
\usepackage{amsfonts}
\usepackage{amscd}
\usepackage{amssymb}
\usepackage{hyperref}
\newcommand{\G}{\Gamma}
\newcommand{\cstar}{\bb{C}^\times}
\newcommand{\cf}{cf.\ }
\newcommand{\diag}[1]{\mbox{diag}\left\{#1\right\}}
\newcommand{\eg}{e.g.}
\newcommand{\lc}{loc.\ cit.}
\newcommand{\eqr}[1]{(\ref{eq:#1})}
\newcommand{\ie}{i.e.\ }

\newcommand{\e}[1]{\textbf{e}\left(\textstyle#1\right)}
\newcommand{\mc}[1]{\mathcal{#1}}
\newcommand{\Teq}{T=\left(\begin{array}{cr}1&1\\0&1\end{array}\right)}
\newcommand{\Seq}{S=\left(\begin{array}{cr}0&-1\\1&0\end{array}\right)}
\newcommand{\abcd}{\begin{pmatrix}a&b\\c&d\end{pmatrix}}
\newcommand{\twomat}[4]{\left(\begin{array}{rr}#1&#2\\#3&#4\end{array}\right)}
\newcommand{\threemat}[9]{\left(\begin{array}{ccc}#1&#2&#3\\#4&#5&#6\\#7&#8&#9\end{array}\right)}

\newcommand{\cvec}[3]{\begin{pmatrix}#1\\#2\\#3\end{pmatrix}}
\newcommand{\bb}[1]{\mathbb{#1}}
\newcommand{\gln}[1]{GL_{#1}(\bb{C})}
\newcommand{\m}{\mc{M}}
\newcommand{\mq}{\m_{\bb{Q}}}
\newcommand{\h}{\bb{H}}
\newcommand{\C}{\bb{C}}
\newcommand{\Q}{\bb{Q}}
\newcommand{\hrho}{\mc{H}(\rho)}
\newcommand{\hrhoq}{\mc{H}(\rho)_\bb{Q}}

\newcommand{\pubd}{$p$\,-unbounded}

\newtheorem{thm}{Theorem}[section]
\newtheorem{lem}[thm]{Lemma}

\newtheorem{cor}[thm]{Corollary}
\newtheorem{prop}[thm]{Proposition}
\newtheorem{conj}[thm]{Conjecture}

\numberwithin{equation}{section}

\begin{document}

\title[three-dimensional Fourier coefficients]{Fourier coefficients of three-dimensional vector-valued modular forms}
\author{Christopher Marks}
\address{Department of Mathematics and Statistics, California State University, Chico}
\email{cmarks@csuchico.edu}

\maketitle
\begin{abstract}
We prove that only a finite number of three-dimensional, irreducible representations of the modular group admit vector-valued modular forms with bounded denominators. This provides a verification, in the three-dimensional setting, of a conjecture concerning the Fourier coefficients of noncongruence modular forms, and reinforces the understanding from mathematical physics that when such a representation arises in rational conformal field theory, its kernel should be a congruence subgroup of the modular group.
\end{abstract}
\tableofcontents
\section{Introduction}
It has been understood, at least since the time of Hecke, that modular forms for congruence subgroups have Fourier coefficients with bounded denominators. In other words, if $G$ is a congruence subgroup of $\G=SL_2(\bb{Z})$ and $f$ is an integral weight cusp form for $G$, whose Fourier expansion at infinity has rational numbers as coefficients, then for some large enough integer $M$ the Fourier coefficients of $Mf$ are integers (see \cite[thm 3.52]{Shimura} for a proof). This is, of course, one of the main reasons that congruence subgroups feature prominently in the theory of modular forms: when the $q$-expansions arising here are integral, these integers represent important quantities in number theory, geometry, or physics.

A physical example of this phenomenon is found in rational conformal field theory (RCFT) and its mathematical  counterpart, the theory of rational vertex operator algebras (VOAs). A rational VOA has associated to it a finite number of irreducible modules, which are $\bb{Z}$-graded complex vector spaces
$$M^{(j)}=\bigoplus_{n\geq0}M_n^{(j)},\ 1\leq j\leq d,$$
whose finite-dimensional summands encode the physical data of the underlying RCFT. Each such module defines a  graded character
$$\chi_j=q^{h_j-\frac{c}{24}}\sum_{n\geq0}\dim M_n^{(j)}q^n,$$
where $h_j\in\bb{Q}$ is the conformal weight associated to $M^{(j)}$ and $c\in\bb{Q}$ is the central charge of the theory. Zhu has shown \cite{Zhu} that if the formal variable $q$ is interpreted as in the theory of modular forms, the $\chi_j$ become holomorphic functions on the complex upper half-plane, and their $\bb{C}$-linear span carries a representation of $\G$. In other words, this space of functions defines a $d$-dimensional vector-valued modular function for $\G$. It has long been conjectured by physicists \cite{Moore} that the kernel of the representation arising in this situation is a congruence subgroup of $\G$, so that the $\chi_j$ are congruence modular functions, and mathematically this is somewhat established \cite{Bantay,Dong}.

Returning to the general situation, it was already understood by Fricke and Klein that this bounded denominator property need not hold when the group of invariance is noncongruence. Indeed, one finds in \cite{Fricke1} the example of 
$$u(\tau)=\int_{i\infty}^\tau\eta^4(z)\,dz=\sum_{n\geq0}\frac{\Psi(n)}{6n+1}q^{n+\frac{1}{6}},$$
where $\eta=q^\frac{1}{24}\prod_{n\geq1}(1-q^n)$ denotes Dedekind's eta function and the integer $\Psi(n)$ counts integral points on the elliptic curve $y^2=-3x^3+n$. In this case the group fixing $u$ is of infinite index in $\G$ (thus noncongruence), and as is indicated by the above expansion there are infinitely many primes $p\equiv1\pmod{6}$ appearing in the denominators of its Fourier coefficients \cite{FM}. When such a form $f$ has rational Fourier coefficients, yet there is no sufficiently large integer $M$ such that the Fourier coefficients of $Mf$ are integral, one says that\footnote{This language is, of course, also used in the more general situation, where $\bb{Q}$ and $\bb{Z}$ are replaced by an arbitrary number field and its ring of integers, respectively. In this article, however, we consider $\bb{Q}$-rational Fourier coefficients only.} $f$ has \emph{unbounded denominators}. More specifically, if there is a prime number $p$ which occurs to an arbitrarily high power in the denominators of the Fourier coefficients of $f$, then we say that $f$ is \emph{\pubd}.

Work of Atkin and Swinnerton-Dyer \cite{ASD} served to rekindle interest in noncongruence modular forms, and subsequent results by Scholl \cite{Scholl1,Scholl2}, Li and Long \cite{Long2,Long1}, and Mason \cite{FM,M2} support the idea that the bounded denominator property completely characterizes modular forms on congruence subgroups. This may be formalized in the following
\begin{conj}\label{conj:ubd}
Suppose $f$ is a modular form for a finite index subgroup $G\leq\G$, which has bounded denominators. Then $f$ is modular for the congruence closure of $G$ (\ie the intersection of all congruence subgroups containing $G$).
\end{conj}
This conjecture extends naturally to the vector-valued setting, and indeed is profitably studied there. The extension formulated by Mason \cite{M2} is
\begin{conj}\label{conj:vvmf} Suppose $F$ is a vector-valued modular form for a representation $\rho$ of $\G$, such that the components of $F$ have Fourier expansions with bounded denominators. Then the kernel of $\rho$ is a congruence subgroup of $\G$.
\end{conj}
From this perspective, the group $G$ of conjecture \ref{conj:ubd} becomes the kernel of the (finite image) representation $\rho$ in conjecture \ref{conj:vvmf} (or $\ker\rho$ is the intersection of the conjugates of $G$ if $G$ is not normal in $\G$), and the components of any vector associated to $\rho$ are (by definition) modular forms for $G$. But of course $\rho$ need not have finite image, in which case the components of any associated vector-valued modular form may be a more general type of function than a traditional ``scalar'' modular form. For example, from the above function $u$ one obtains a vector $(u,1)^t$ which transforms according to a representation $\rho:\G\rightarrow\gln{2}$ whose kernel (the aforementioned fixing group for $u$) is the normal closure of $\scriptsize{\twomat{1}{6}{0}{1}}$ in $\G$, which is called $\Delta(6)$ by Wohlfahrt \cite{W} and is of infinite index in $\G$. 

In light of the modular invariance of graded characters discussed above, a strong motivation for proving conjecture \ref{conj:vvmf} is that it would provide a mathematical verification, using only intrinsic properties of modular forms, of the expectation from RCFT that the representations of $\G$ arising here have  congruence kernels. There is no doubt that a proof of the above conjectures will necessarily involve some extremely subtle, and arithmetically rich, ideas from the theory of Riemann surfaces. It is reasonable to expect that these ideas would inform physics as well as number theory, in ways that cannot be fully predicted.

Conjecture \ref{conj:vvmf} is easily verified when the representation is one-dimensional, for here, as is well known, there are only 12 possible representations to consider, each of whose kernel is congruence (of level dividing 12). And in the two-dimensional setting, Mason \cite{M2} has verified conjecture \ref{conj:vvmf} for all but a finite number of open cases (which have recently been settled \cite{FM}). In this article we shall, in a similiar fashion, verify conjecture \ref{conj:vvmf} in dimension three. Explicitly, we prove
\begin{thm}\label{thm:main}
Up to equivalence of representation, only a finite number of irreducible $\rho:\G\rightarrow\gln{3}$, with $\rho\scriptsize{\twomat{1}{1}{0}{1}}$ of finite order, admit vector-valued modular forms with bounded denominators.
\end{thm}
 
An outline of the remainder of this article is as follows. In the next section, we review the necessary topics from the theory of vector-valued modular forms. In section \ref{sec:3d} we utilize a recursive formula from the Fuchsian theory of ordinary differential equations, to establish that the generators for spaces of three-dimensional vector-valued modular forms almost always have unbounded denominators. This allows us to complete, in section \ref{sec:main}, the proof of theorem \ref{thm:main}. In section \ref{sec:monomial}, we apply thereom \ref{thm:main} to the classical setting, where the representation $\rho$ has finite image. This application supports conjecture \ref{conj:ubd}, by exhibiting infinitely many new families of noncongruence modular forms with unbounded denominators.  

\section{Background}\label{sec:background}
In this section, we recall material from the theory of vector-valued modular forms and Fuchsian differential equations needed to establish theorem \ref{thm:main}. For more details regarding the theory of vector-valued modular forms, in addition to the references listed in this section one may consult the author's doctoral dissertation \cite{Marks1}. The facts cited below pertaining to the theory of Fuchsian differential equations are included in any elementary text on the subject, \eg\ the reference \cite{Hille} given below.

Let $\h$ denote the complex upper half-plane, $\mc{H}$ the complex linear space of holomorphic functions
$f:\h\rightarrow\C$, and $\G=SL_2(\bb{Z})$ the full modular group of $2\times2$ matrices with integer entries and determinant 1. We denote by $\scriptsize{\Seq}$, $\scriptsize{\Teq}$ the well-known generators of $\G$. For each integer $k$ we write
\begin{eqnarray*}
|_k:\mc{H}\times\G&\rightarrow&\mc{H},\\
(f,\gamma)&\mapsto&f|_k\gamma
\end{eqnarray*}
to denote the $k^{th}$ slash action of $\G$ on $\mc{H}$; thus for each $\gamma=\scriptsize{\abcd}\in\G$, $\tau\in\h$ we have
$$f|_k\gamma(\tau)=(c\tau+d)^{-k}f\left(\frac{a\tau+b}{c\tau+d}\right).$$
A holomorphic function $F:\h\rightarrow\bb{C}^m$ is a $m$-dimensional \emph{vector-valued modular form} of weight $k\in\bb{Z}$ if
the component functions comprising $F$ satisfy a moderate growth condition at the cusps of $\G$, just as in the
classical theory of modular forms, and if the span of these components is an invariant subspace of $\mc{H}$ under the $|_k$ action
of $\G$ on $\mc{H}$. Explicitly, if one writes $F$ as a column vector $F=(f_1,\cdots,f_m)^t$ then the above action of $\G$ on the span of the $f_j$ takes the form of a matrix representation $\rho:\G\rightarrow\gln{m}$, and we say that $F$ is a vector-valued modular form of weight $k$ for $\rho$ if the following conditions are satisfied:
\begin{enumerate}\label{def:vvmf}
\item The components $f_j$ of $F$ are of \emph{moderate growth at infinity}, \ie there is an integer $N$ such that for each $j$
we have $|f_j(x+iy)|<y^N$ for any fixed $x\in\bb{R}$ and $y\gg0$.
\item The functional equation $F|_k\gamma=\rho(\gamma)F$ is satisfied for each $\gamma\in\G$ (here $|_k$ is applied componentwise to $F$).
\end{enumerate}
We write $\mc{H}(k,\rho)$ for the complex linear space of all such vectors. Regardless of $\rho$, if $k$ is large enough then $\mc{H}(k,\rho)\neq0$ \cite[cor 3.12]{KM4}, and there is a minimal weight $k_0\geq1-m$ such that $\mc{H}(k_0,\rho)\neq0$.  Furthermore, if $\rho$ is indecomposable -- as shall be the case throughout this article -- then it follows directly from the above definition that $\rho(S^2)=(-1)^k$ whenever $\mc{H}(k,\rho)\neq0$. Thus in the indecomposable setting there is a $\bb{Z}$-graded space
\begin{equation}\label{eq:hrho}
\mc{H}(\rho)=\bigoplus_{k\geq0}\mc{H}(k_0+2k,\rho)
\end{equation}
containing all holomorphic, integral weight vector-valued modular forms for $\rho$. If $\rho=\textbf{1}$ is the trivial one-dimensional representation of $\G$, then\footnote{More generally, if $\rho$ is unitary then the minimal weight $k_0$ is nonnegative.} $k_0=0$ and we write
$$\mc{H}(\textbf{1})=\m=\bigoplus_{k\geq0}\m_{2k}$$
for the graded ring of holomorphic, integral weight modular forms for $\G$. As is well-known, $\m=\C[E_4,E_6]$ is a graded polynomial algebra in $E_4\in\m_4$ and $E_6\in\m_6$, where for each even integer $k\geq2$ we write
\begin{equation}\label{eq:Eisenstein}
E_k(\tau)=1-\frac{2k}{B_k}\sum_{n\geq1}\sigma_{k-1}(n)q^n
\end{equation}
for the normalized Eisenstein series in weight $k$; here $q=q(\tau)=e^{2\pi i\tau}$, $B_k$ denotes the $k^{th}$ Bernoulli number, and $\sigma_k(n)=\sum_{1\leq d|n}d^k$. Each space \eqr{hrho} is a graded $\m$-module via componentwise multiplication. 

If $\rho'$ is equivalent to $\rho$ -- meaning there is a $U\in\gln{m}$ such that $\rho'(\gamma)=U\rho(\gamma)U^{-1}$ for each $\gamma\in\G$ -- then multiplication by $U$ defines an isomorphism $\mc{H}(\rho)\cong\mc{H}(\rho')$ of graded $\m$-modules. This allows one to study vector-valued modular forms for representations having desirable matrix forms within their particular equivalence class. In particular, in this article we shall be concerned only with those $\rho$ such that $\rho(T)$ is of finite order. Thus the above isomorphism allows us to assume that
\begin{equation}\label{eq:rhoT}
\rho(T)=\diag{\e{r_1},\cdots,\e{r_m}}
\end{equation}
for some rational numbers $0\leq r_j<1$; here and throughout we write $\e{r}=e^{2\pi ir}$ for the exponential of $r\in\bb{R}$. In this case, the moderate growth condition implies that each $F\in\mc{H}(\rho)$ has a holomorphic $q$-expansion
\begin{equation}\label{eq:F}
F(\tau)=\cvec{f_1(\tau)}{\vdots}{f_m(\tau)}=\cvec{\sum_{n\geq0}a_1(n)q^{r_1+n}}{\vdots}{\sum_{n\geq0}a_m(n)q^{r_m+n}},
\end{equation}
with $a_j(n)\in\bb{C}$ for each $j,n$. If the $a_j(n)$ are rational numbers, then we call $F$ \emph{\pubd} whenever at least one component of $F$ is.

A fundamental fact concerning the $\m$-module structure of \eqr{hrho} is 
\begin{thm}\label{thm:free}
If $\rho$ is indecomposable and $\rho(T)$ has finite order, then $\mc{H}(\rho)$ is a free $\m$-module of rank equal to the dimension of $\rho$.\qed
\end{thm}
\noindent See \cite[Thm 1]{MM} for a proof, or \cite{KM4} for a more general result.

Recall that the modular derivative in weight $k\in\bb{Z}$ is the operator
\begin{eqnarray}\label{eq:Dk}
D_k&:&\mc{H}\rightarrow\mc{H},\nonumber\\
D_kf&=&\frac{1}{2\pi i}\frac{df}{d\tau}-\frac{k}{12}E_2f,
\end{eqnarray}
with $E_2$ as in \eqr{Eisenstein}. It is well-known \cite[sec 10.5]{Lang} that this derivative is covariant with respect to the slash action of the modular group, so that 
$$(D_kf)|_{k+2}\gamma=D_k(f|_k\gamma)$$
for any meromorphic $f:\bb{H}\rightarrow\bb{C}$, $\gamma\in\G$, and $k\in\bb{Z}$. This covariance helps establish the important fact that $D_k$ takes (quasi-)modular forms of weight $k$ to those of weight $k+2$, and in particular we have 
\begin{equation}\label{eq:E246}
D_2E_2=-\frac{1}{12}[E_2^2-E_4],\ \ D_4E_4=-\frac{1}{3}E_6,\ \ D_6E_6=-\frac{1}{2}E_4^2.
\end{equation}
This generalizes to higher dimension and yields a weight two
operator $D$, which acts (componentwise) on each graded space \eqr{hrho} of vector-valued modular forms by acting as $D_k$ on $\mc{H}(k,\rho)$. One defines for any $n\geq2$ the composition
\begin{equation}\label{eq:Dkn}
D_k^n=D_{k+2(n-1)}\circ\cdots\circ D_k
\end{equation}
and in this way powers of $D$ are well-defined operators on \eqr{hrho}. This allows one to define a skew polynomial ring $\mc{R}$, which as an additive group is just the polynomial ring in one variable $\m[d]$, and whose multiplication is defined by the identity $dM=Md+DM$ for each $M\in\m$. Each space \eqr{hrho} of vector-valued modular forms is then a graded left $\mc{R}$-module (finitely generated thanks to Theorem \ref{thm:free}), where $\m$ again acts by componentwise multiplication, and $d^n$ acts as the $n^{th}$ power of the modular derivative, \ie $d^n$ acts as \eqr{Dkn} on $\mc{H}(k,\rho)$. 

A \emph{modular differential equation} is simply an equation $L[f]=0$ with $L\in\mc{R}$ homogeneous by weight. A very special case -- which includes the setting of the present article -- occurs when $L$ is monic, so that
\begin{equation}\label{eq:monic}
L[f]=D_k^mf+M_4D_k^{m-2}f+\cdots+M_{2m}f=0
\end{equation}
for some $m\geq1$, $k\in\bb{Z}$, and $M_j\in\m_j$ for each $j$. From the covariance of the modular derivative, it follows that the space of solutions of \eqr{monic} is invariant under the $|_k$ action of $\G$, thus such equations yield candidates for vector-valued modular forms. Because $E_2$ and the $M_j$ are holomorphic in $\h$, one sees from \eqr{E246} and \eqr{Dkn} that the only singular point of \eqr{monic} is $i\infty$, \ie $q=0$. Using the change of variable $\frac{1}{2\pi i}\frac{d}{d\tau}=q\frac{d}{dq}$, one obtains from \eqr{monic} an equation
\begin{equation}\label{eq:qversion}
q^m\frac{df}{dq^m}+g_{m-1}(q)q^{m-1}\frac{df}{dq^{m-1}}+\cdots+g_0(q)f=0
\end{equation}
for some functions 
\begin{equation}\label{eq:gj}
g_j(q)=\sum_{n\geq0}G_j(n)q^n\in\bb{C}[E_2,E_4,E_6]
\end{equation}
which are holomorphic in the disk $|q|<1$. One sees from \eqr{qversion} that in fact $q=0$ is a regular singular point in the sense of Fuchs \cite[Chs 5,9]{Hille}, so  that solutions of \eqr{qversion} may be obtained by employing the well-known recursive formula of Fuchs and Frobenius. To implement this method, a normalized series solution of the form
\begin{equation}\label{eq:formal}
f(q)=q^r\left[1+\sum_{n\geq1}a(n)q^n\right]
\end{equation}
is assumed and evaluated according to \eqr{qversion}. Since the resulting function $L[f]$ must be identically zero, one obtains an infinite set of conditions which must be satisfied by $r$ and the $a(n)$. The first condition is known classically as the \emph{indicial equation} $\phi(r)=0$ with 
\begin{eqnarray}\label{eq:indicial}
\phi(r)&=&r(r-1)\cdots(r-(m-1))\nonumber\\
&\ &+\,G_{m-1}(0)r(r-1)\cdots(r-(m-2))\\
&\ &+\cdots+G_1(0)r+G_0(0),\nonumber
\end{eqnarray}
and it determines the leading exponents of the solutions of \eqr{qversion} uniquely; these values of $r$ are called the \emph{indicial roots} of \eqr{qversion}. Continuing on, one finds at the $n^{th}$ step that the condition
\begin{equation}\label{eq:recursivegen}
a(n)=\frac{a(n-1)\phi_1(r+(n-1))+\cdots+a(1)\phi_{n-1}(r+1)+\phi_n(r)}{\phi(r+n)}
\end{equation}
must be satisfied, where
\begin{eqnarray*}
\phi_j(r)&=&G_{m-1}(j)r(r-1)\cdots(r-(m-2))\\
&\ &\,+G_{m-2}(j)r(r-1)\cdots(r-(m-3))\\
&\ &+\,\cdots+G_1(j)r+G_0(j).
\end{eqnarray*}
Thus \eqr{recursivegen} defines each coefficient of \eqr{formal} recursively, so long as there does not exist an $n\geq1$ such that $P(r+n)=0$. In other words, so long as no two indicial roots of \eqr{qversion} differ by an integer one derives from this method a basis $\{f_1,\cdots,f_m\}$ of formal solutions of \eqr{qversion}, each of the form \eqr{formal}, and \cite[sec 9.1]{Hille} for each $j$ the ratio $\frac{f_j(q)}{q^{r_j}}$ converges to a holomorphic function in the open disk $|q|<1$. 

The upshot of the above discussion is that if the indicial roots $r_1,\cdots,r_m$ of \eqr{qversion} are nonnegative real numbers which are distinct$\pmod{\bb{Z}}$, then there is a representation $\rho:\G\rightarrow\gln{m}$ and a vector-valued modular form $F\in\mc{H}(k,\rho)$ of the form \eqr{F}, with $a_j(0)=1$ for each $j$, such that the $f_j$ form a fundamental system of solutions (\ie a basis of the solution space) of \eqr{qversion}. Note that this choice of basis implies that $\rho(T)$ is diagonal as in \eqr{rhoT}, so that the indicial roots of \eqr{monic} may be taken as the exponents of the eigenvalues of $\rho(T)$.

It is of fundamental importance for the analysis undertaken in this article that the space \eqr{hrho} associated to any irreducible $\rho:\G\rightarrow\gln{3}$ with $\rho(T)$ as in \eqr{rhoT} is a cyclic $\mc{R}$-module \cite[thm 3]{MM}, whose generator $F$ has components spanning the solution space of a monic differential equation \eqr{monic}. Furthermore, if $\rho(T)$ has finite order (\ie the indicial roots of \eqr{monic} are rational numbers), then the Fourier coefficients of $F$ will be rational numbers as well. In this setting one may exploit the recursive formula \eqr{recursive} to find sufficient conditions on the indicial roots so that $F$ will have unbounded denominators. Because the generators of \eqr{hrho} as $\m$-module are $F,DF,D^2F$, and using the fact that the Eisenstein series \eqr{Eisenstein} have bounded denominators, one may deduce the unboundedness of the denominators of an arbitrary vector-valued modular form in \eqr{hrho} from that of $F$. This is the strategy employed below. 

Following Wohlfahrt \cite{W}, we refer to the order of $\rho(T)$ in $\rho(\G)$ as the \emph{level} of $\rho$. If $\rho$ is of finite level $N$ then $\ker\rho$ is a normal subgroup of $\G$ that contains $T^N$, thus it also contains the normal closure $\Delta(N)$ of the subgroup of $\G$ generated by $T^N$. It is proven in \lc\ that if $N<6$ then $\Delta(N)$ is the principal congruence subgroup $\G(N)$, so $\ker\rho$ is necessarily a congruence subgroup of level $N$ in this case; in particular, the components of any vector-valued modular form for a representation of level less than six are congruence modular forms, and consequently have bounded denominators. On the other hand, it is also proven in \lc\ that $\Delta(N)$ is of infinite index in $\G(N)$ when $N\geq6$, so in this setting the image of $\rho$ may be finite or infinite.

\section{Three-dimensional vector-valued modular forms}\label{sec:3d}
In this section we analyze the Fourier coefficients of the minimal weight vector-valued modular form \eqr{F0} associated to a generic three-dimensional representation of the modular group. 
This analysis forms the core of the proof of theorem \ref{thm:main}, which will be completed in section \ref{sec:main}.

Suppose that
\begin{equation}\label{eq:rho3d}
\rho:\G\rightarrow\gln{3}
\end{equation}
is irreducible such that $\rho(T)$ has finite order. Up to equivalence of representation, we may and shall now assume that $\rho(T)$ is diagonal as in \eqr{rhoT}. It follows directly from \cite[thm 1]{MM} (or see \cite[thm 4.1]{Marks2} for a proof) that the space \eqr{hrho} of holomorphic vector-valued modular forms for $\rho$ is a cyclic $\mc{R}$-module
\begin{equation}\label{eq:free3}
\hrho=\mc{R}F_0=\m F_0\oplus\m DF_0\oplus\m D^2F_0,
\end{equation}
and by \cite[thm 3]{MM} the components of the generator $F_0$ form a fundamental system of solutions of a modular differential equation
\begin{equation}\label{eq:mde3}
L[f]=D_{k_0}^3f+\alpha_4E_4D_{k_0}f+\alpha_6E_6f=0.
\end{equation}
Here $k_0=4r-2\in\bb{Z}$, with $r=r_1+r_2+r_3$ the sum of the exponents in \eqr{rhoT} (which are also the indicial roots of \eqr{mde3}), $E_k$ is the Eisenstein series \eqr{Eisenstein}, and the complex numbers $\alpha_4,\alpha_6$ are uniquely determined by the $r_j$; \cf \cite[lem 2.3]{Marks2}, or see \eqr{walpha} for the explicit formulae. Note also that the $r_j$ are distinct by \cite[thm 3]{MM}. Denoting the order of $\rho(T)$ by $N$, we have
\begin{equation}\label{eq:rhoT3}
\rho(T)=\diag{\e{\frac{A}{N}},\e{\frac{B}{N}},\e{\frac{C}{N}}},
\end{equation}
where the integers $A,B,C$ are distinct and satisfy
\begin{equation}\label{eq:ABC}
0\leq A,B,C\leq N-1,\ \ (A,B,C,N)=1.
\end{equation}
Knowing this, we may now assume that the minimal weight vector $F_0\in\mc{H}(k_0,\rho)$ has a Fourier expansion of the form
\begin{equation}\label{eq:F0}
F_0=\cvec{q^\frac{A}{N}+\sum_{n\geq1}a(n)q^{\frac{A}{N}+n}}{q^\frac{B}{N}+\sum_{n\geq1}b(n)q^{\frac{B}{N}+n}}{q^\frac{C}{N}+
\sum_{n\geq1}c(n)q^{\frac{C}{N}+n}}.
\end{equation}
Using \eqr{Dk},\eqr{E246}, and \eqr{Dkn}, one obtains from \eqr{mde3} an equation \eqr{qversion} with $m=3$ and 
\begin{eqnarray}\label{eq:g3}
g_2(q)&=&\sum_{n\geq0}G_2(n)q^n\nonumber\\\nonumber\\
&=&3+(3k_0+6)P,\nonumber\\ \nonumber\\
g_1(q)&=&\sum_{n\geq0}G_1(n)q^n\\
&=&1+(3k_0+6)P+(3k_0^2+9k_0+6)P^2\nonumber\\
&\ &+(3k_0+2+144\alpha_4)Q,\nonumber\\ \nonumber\\
g_0(q)&=&\sum_{n\geq0}G_0(n)q^n\nonumber\\
&=&k_0(3k_0+2+144\alpha_4)PQ\nonumber+k_0(k_0+1)(k_0+2)P^3\nonumber\\
&\ &+(k_0-432\alpha_6)R\nonumber;
\end{eqnarray}
here we set
\begin{equation}\label{eq:pqr}
P=-\frac{1}{12}E_2,\ \ Q=\frac{1}{144}E_4,\ \ R=-\frac{1}{432}E_6.
\end{equation}
Using this notation (and setting $a(0)=1$), we may write the Fuchsian recursive relation \eqr{recursivegen} for the coefficients of the first component of \eqr{F0} as
\begin{equation}\label{eq:recursive}
a(n)=-\frac{1}{\phi\left(\frac{A}{N}+n\right)}\sum_{j=0}^{n-1}a(j)\phi_{n-j}\left(\frac{A}{N}+j\right),
\end{equation}
where
\begin{eqnarray}\label{eq:indicial}
\phi_j(\lambda)&=&G_2(j)\lambda(\lambda-1)+G_1(j)\lambda+G_0(j),\\
\phi(\lambda)&=&\lambda(\lambda-1)(\lambda-2)+\phi_0(\lambda).\nonumber
\end{eqnarray}
The exponents $\frac{A}{N},\frac{B}{N},\frac{C}{N}$ in \eqr{rhoT3} are the solutions of the indicial equation $\phi(\lambda)=0$ associated to \eqr{mde3}, and from this one obtains directly
\begin{eqnarray*}
G_2(0)&=&3-\frac{\sigma}{N},\\ \\
G_1(0)&=&G_2(0)+\frac{\omega}{N^2}-2,\\ \\
G_0(0)&=&-\frac{\varpi}{N^3},
\end{eqnarray*}
where we set
$$\sigma=A+B+C,\ \ \omega=AB+AC+BC,\ \ \varpi=ABC.$$
Using this information, it is now straightforward to compute and find
\begin{equation}\label{eq:phi}
\phi\left(\frac{A}{N}+n\right)=\frac{n\lambda(n)}{N^2},
\end{equation}
where for each $n\geq1$ we set
\begin{equation}\label{eq:lambdan}
\lambda(n)=Nn[Nn+(A-B)+(A-C)]+(A-B)(A-C).
\end{equation}
Furthermore, comparing the above formulae for the $G_j(0)$ with those obtained directly from \eqr{g3} yields
\begin{eqnarray}\label{eq:walpha}
k_0&=&\frac{x_0}{N},\nonumber\\
\nonumber\\
\alpha_4&=&\frac{x_4}{(12N)^2},\\
\nonumber\\
\alpha_6&=&\frac{x_6}{(12N)^3},\nonumber
\end{eqnarray}
where we define the integers
\begin{eqnarray}\label{eq:xj}
x_0&=&4\sigma-2N,\nonumber\\
\nonumber\\
x_4&=&144\omega+x_0(12N-3x_0)+8N^2,\\
\nonumber\\
x_6&=&x_0x_4+x_0(x_0+2N)(x_0+4N)-1728\varpi.\nonumber
\end{eqnarray}
A final round of elementary computations yields
\begin{eqnarray*}
G_2(1)&=&\frac{24\sigma}{N},\\ \\
G_1(1)&=&\frac{240\omega-48\sigma(2\sigma-N)}{N^2},\\ \\
G_0(1)&=&\frac{504\varpi+\big(2\sigma-N\big)\big(8\sigma(4\sigma-N)-120\omega\big)}{N^3},\\
\end{eqnarray*}
and from this it is trivial to verify the following
\begin{lem}\label{lem:zn}
For each $n\geq0$ we have $\phi_1\left(\frac{A}{N}+n\right)=\frac{z_n}{N^3}$, where
\begin{eqnarray}\label{eq:zn}
z_n&=&24[10\omega Nn+\sigma(A+Nn)(A+N(n-1))]+\nonumber\\
& &8[2\sigma-N][\sigma(4\sigma-N)-15\omega-6\sigma(A+Nn)]+\\
& &240A\omega+504\varpi.\nonumber
\end{eqnarray}
\end{lem}\qed

It is easy to see that the $x_j$ in \eqr{xj} satisfy $2|x_0$, $4|x_4$, $8|x_6$, and this will be used to prove 
\begin{lem}\label{lem:Gj}
For $j=2,3$ set
$$\delta_j=\left\{\begin{array}{lr}0&j\mid N,\\1&j\nmid N.\end{array}\right.$$
Then for each $n\geq2$ the following hold:
\begin{enumerate}
\item $G_2(n)\in\bb{Z}$.\\
\item $G_1(n)\in\frac{1}{3^{\delta_3}N^2}\bb{Z}$.\\
\item $G_0(n)\in\frac{1}{2^{\delta_2}3^{\delta_3}N^3}\bb{Z}$.\\
\end{enumerate}
\end{lem}

\proof Since the minimal weight $k_0$ is an integer, one sees from \eqr{walpha}, \eqr{xj} that $N|4\sigma$. Furthermore, from \eqr{g3} and \eqr{Eisenstein} we obtain
$$G_2(n)=\frac{24\sigma}{N}\sigma_1(n)$$
for each $n\geq1$. This implies statement $(1)$. For statements $(2)$ and $(3)$, it is sufficient to verify the analogous statement for the coefficients of each summand of $g_1$ and $g_2$, respectively. This amounts to a routine verification and we omit the proof. We do note, however, that in addition to the 2-adic properties of the $x_j$ mentioned above, we also have that $3|x_4$ iff $3|N$, and $3|N$ implies $3|x_0$, $9|x_6$; these observations are all that is required to fill in the remaining details. \qed

For each prime $p$ we write $\nu_p$ to denote the $p$-adic valuation of $\bb{Q}$; thus if $x=p^kv$ is an integer with $(v,p)=1$ then $\nu_p(x)=k$, and if $\frac{x}{y}\in\bb{Q}$ then $\nu_p\left(\frac{x}{y}\right)=\nu_p(x)-\nu_p(y)$. The most important step in the proof of theorem \ref{thm:main} is the determination of $\nu_p(z_n)$ for various primes $p$, and we turn now to this task.
\begin{prop}\label{prop:xn}
Let $p$ be a prime dividing the level $N$ of \eqr{rho3d}. After relabeling (if needed) the indicial roots $\frac{A}{N},\frac{B}{N},
\frac{C}{N}$ of \eqr{mde3}, the following statements hold for all $n\geq0$:
\begin{enumerate}
\item If $p>7$, then $\nu_p(z_n)=0$.
\item If $p=7$ and $7\mid\varpi$, then $\nu_7(z_n)=0$.
\item If $p=7$ and $7\nmid\varpi$, $7^2|N$, then one of the following holds:
\begin{enumerate}
\item $7\mid\omega$ and $\nu_7(z_n)=1$.
\item $7\nmid\omega$ and $\nu_7(z_n)=0$.
\end{enumerate}
\item If $p=5$ and $5\nmid\varpi$, then $\nu_5(z_n)=0$.
\item If $p=5$ and $5|\varpi$, $5^2|N$, then $\nu_5(z_n)=1$.
\item If $p=3$ and $3^2|N$, $3\nmid\omega$, then $\nu_3(z_n)=1$.
\item If $p=3$ and $3^3|N$, $3|\omega$, then $\nu_3(z_n)=2$.
\item If $p=2$ and $2^5|N$, then $\nu_2(z_n)=4$.
\end{enumerate}
\end{prop}

\proof From the formula \eqr{walpha} for $k_0$, it follows that if $p\geq3$ is a prime dividing $N$ then $\nu_p(\sigma)\geq\nu_p(N)>0$,
and if $p=2$ then $\nu_2(\sigma)\geq\nu_p(N)-2$. On the other hand, it is seen from \eqr{zn} that if $p$ divides both $N$ and $\sigma$, then $\nu_p(z_n)>0$ iff $p$ divides
\begin{equation}\label{eq:term}
240A\omega+504\varpi=24A[10A(B+C)+31BC].
\end{equation}
By assuming, as we may, that $p\nmid AB$, it is clear that if $p|C$ and $p>5$, then $p$ does not divide \eqr{term}; this
implies statement $(2)$ and part of statement $(1)$ of the proposition.

Assume now that $\nu_p(N)=k\geq1$ with $5\leq p\nmid\varpi$. Then \eqr{term} shows that $\nu_p(z_n)>0$ iff $p|10A(B+C)+31BC$.
Transposing $A$ and $B$ throughout the calculations which led to \eqr{zn} will yield the analogue $y_n$ of \eqr{zn} for the numerator of $\phi_1\left(\frac{B}{N}+n\right)$, and one finds similarly that $\nu_p(y_n)>0$ iff $p|10B(A+C)+31AC$. Noting the fact we also have
$\nu_p(\sigma)\geq k$ in this context, a trivial calculation shows that for any $1\leq m\leq k$ we have
\begin{eqnarray*}
10A(B+C)+31BC&\equiv&-(10A^2+31AB+31B^2)\pmod{p^m},\\
10B(A+C)+31AC&\equiv&-(10B^2+31AB+31A^2)\pmod{p^m},
\end{eqnarray*}
so a necessary condition for $\nu_p(z_n)\geq m$ and $\nu_p(y_n)\geq m$ to \emph{both} hold is that $p^m$ divides the difference
\begin{equation}\label{eq:diffsq}
10A^2+31AB+31B^2-(10B^2+31AB+31A^2)=21(B^2-A^2).
\end{equation}
Note that $p\nmid(B+A)$, since $p|\sigma$ and $p\nmid C$, so $p^m|(B^2-A^2)$ iff $A\equiv B\pmod{p^m}$.
Furthermore, if this holds then it follows immediately that
\begin{equation}\label{eq:termmodp}
10A(B+C)+31BC\equiv-72A^2\pmod{p^m},
\end{equation}
and because we are assuming $5\leq p\nmid A$, this cannot be. In particular, this implies that $p$ does not divide
\eqr{termmodp} if $p=5$ or $p>7$. Taking $m=1$ and relabeling (if needed) then completes the proof of statement $(1)$ of
the Proposition and yields statement $(4)$ as well. On the other hand, assuming $p=7$, $k\geq2$, $m=2$ makes it clear that $\nu_7(z_n)\leq1$ (after relabeling if needed), and from this and \eqr{term} statement $(3)$ follows immediately.

Next we assume $5^2|N$, $5|C$, say $C=5X$ for some integer $X$. Then it follows directly from \eqr{term} that $\nu_5(z_n)\geq1$,
and $\nu_5(z_n)\geq2$ iff $5|(2A+X)$. As in the previous paragraph, we pursue an identical analysis for the integer $y_n$ which
is the numerator of $\phi_1\left(\frac{B}{N}+n\right)$, and find this time that $5^2|y_n$ iff $5|(2B+X)$. Thus $5^2$ necessarily divides the
difference $2(A-B)$ of these terms if $5^2$ divides the numerators of both $a(1)$ and $b(1)$, which is to say $5|(A-B)$. But $5|C$,
$5|\sigma$ imply that $5|(A+B)$, thus $5\nmid(A-B)$ since $(A,B,C)=1$ and $5|C$. This proves Statement $(5)$ of the proposition.

Now assume that $\nu_3(N)=k\geq1$, $3\nmid A$. Statement $(6)$ of the Proposition follows immediately from \eqr{term} by assuming $k\geq2$. On the other hand, if $k\geq3$ and $3|\omega$, then \eqr{term} makes it clear that $\nu_3(z_n)\geq2$.
But examining the calculations which led to \eqr{termmodp}, one sees that this logic remains valid for the prime 3, and taking $m=3$
shows that, up to relabeling, we have $\nu_3(z_n)\leq2$, and this implies statement $(7)$ of the proposition.

Finally, assume that $\nu_2(N)\geq4$. Then $\nu_2(\sigma)\geq2$, and this implies $\nu_2(\varpi)\geq1$. Note that $2\nmid\omega$ since $2|\varpi$. If $\nu_2(\varpi)\geq2$, then we may assume that $2\nmid A$, $2^2|BC$, and this makes it clear that the first two terms of \eqr{zn} are divisible by $2^5$, whereas the last term is divisible only by $2^4$. On the other hand, if $\nu_2(\varpi)=1$, then we may assume that $\nu_2(A)=1$, $2\nmid BC$, and in this case we find that $2^6$ divides the first two terms of \eqr{zn}, but the last term is divisible by only $2^4$. Thus statement $(7)$ holds, and this completes the proof of the proposition.\qed

It follows immediately that if $p$ is a prime satisfying one of conditions $(1)-(8)$ in proposition \ref{prop:xn} then the difference
\begin{equation}\label{eq:delta}
\delta=\delta(p)=\nu_p(z_n)-\nu_p(N)<0
\end{equation}
is well-defined, independently of the integer $n\geq0$. With this notation, we may now prove
\begin{prop}\label{prop:formula}
Suppose $p$ is a prime satisfying one of conditions $(1)-(8)$ of proposition \ref{prop:xn}, and assume furthermore that $\nu_p(N)>2\nu_p(z_0)$. Then for all $n\geq1$ we have
\begin{equation}\label{eq:ord}
\nu_p(a(n))=n\delta-\nu_p\left(\prod_{k=1}^nk\lambda(k)\right),
\end{equation}
with $\lambda(k)$ as in \eqr{lambdan}. In particular, $\nu_p(a(n))$ is a strictly decreasing, negative function of $n$, and \eqr{F0} is \pubd.
\end{prop}

\proof The proof will be made by induction on $n$. Since
$$a(1)=-\frac{N^2}{\lambda(1)}\cdot\phi_1\left(\frac{A}{N}\right)=-\frac{z_0}{N\lambda(1)}$$
with $z_n$ as in Lemma \ref{lem:zn}, it is clear from proposition \ref{prop:xn} that if $p$ is a prime satisfying the hypothesis of the current proposition, then \eqr{ord} holds for $n=1$. Now assume that $n\geq2$ and \eqr{ord} holds for all $1\leq j\leq n-1$. Then   \eqr{recursive}, \eqr{phi}, and basic properties of $\nu_p$ imply that
\begin{eqnarray*}
\nu_p(a(n))&=&\nu_p(N^2)-\nu_p(n\lambda(n))+\nu_p\left(a(n-1)\phi_1\left(\frac{A}{N}+(n-1)\right)\right)\\
&=&n\delta-\nu_p\left(\prod_{k=1}^nk\lambda(k)\right)
\end{eqnarray*}
and the proposition is proved, so long as we have
\begin{equation}\label{eq:ineq}
\nu_p\left(a(n-1)\phi_1\left(\frac{A}{N}+(n-1)\right)\right)<\nu_p\left(a(j)\phi_{n-j}\left(\frac{A}{N}+j\right)\right)
\end{equation}
for all $0\leq j\leq n-2$. Now lemma \ref{lem:Gj} and \eqr{indicial} imply that for each such $j$, there is a $y_j\in\bb{Z}$ such that $\phi_{n-j}\left(\frac{A}{N}+j\right)=\frac{y_j}{2^{\delta_2}3^{\delta_3}N^3}$, and we note that the definition of the $\delta_j$ in lemma \ref{lem:Gj} implies that $\nu_p(2^{\delta_2}3^{\delta_3})=0$ for any $p$ dividing $N$. By the induction hypothesis, it is then sufficient to prove that 
$$(n-(j+1))\delta+\nu_p(z_n)=(n-j)\nu_p(z_n)-(n-(j+1))\nu_p(N)<0$$
for $0\leq j\leq n-2$. It is now apparent that the additional assumption  $\nu_p(N)>2\nu_p(z_n)$ in the statement of the proposition is enough to ensure that \eqr{ineq} holds for each $n$, and this completes the proof.\qed

\begin{cor}\label{cor:ubd}
Suppose there is a prime $p$ which divides $\frac{N}{(N,2^8\cdot3^4\cdot5^2\cdot7^2)}$. Then \eqr{F0} is \pubd.

\proof This amounts to checking the statement of proposition \ref{prop:xn} to see that any such prime $p$ also satisfies the condition $\nu_p(N)>2\nu_p(z_n)$ in proposition \ref{prop:formula}. \qed
\end{cor}

With this corollary in hand, we are now well-situated to complete the proof of theorem \ref{thm:main}.

\section{Proof of main theorem}\label{sec:main}
We now complete the proof of theorem \ref{thm:main}, by establishing some general facts about cyclic $\mc{R}$-modules of vector-valued modular forms.

Assume that $\rho:\G\rightarrow\gln{d}$ is an irreducible representation of arbitrary dimension $d$, with $\rho(T)$ as in \eqr{rhoT}, such that the graded space
$$\hrho=\bigoplus_{j=0}^{d-1}\m D^j F_0=\mc{R}F_0$$
of holomorphic vector-valued modular forms for $\rho$ is a cyclic $\mc{R}$-module with generator $F_0$. We set
\begin{eqnarray*}
\hrhoq&=&\{F\in\hrho\mid F\mbox{ has rational Fourier coefficients}\},\\
\mq&=&\{f\in\m\mid f\mbox{ has rational Fourier coefficients}\}.
\end{eqnarray*}
Then $\hrhoq$ is clearly an $\mq$-module, and we have
\begin{lem}\label{lem:qform}
If $F_0\in\hrhoq$, then
$$\hrhoq=\bigoplus_{j=0}^{d-1}\mq D^jF_0$$
is a free $\mq$-module of rank $d$.
\end{lem}

\proof It follows directly from \eqr{Eisenstein} and \eqr{Dk} that $D^jF_0\in\hrhoq$ for any integer $j\geq0$, so clearly the free $\mq$-module
$\bigoplus_{j=0}^{d-1}\mq D^jF_0$ is contained in $\hrhoq$. On the other hand, suppose
$$F=\cvec{f_1}{\vdots}{f_d}=\cvec{\sum_{n\geq0}c_1(n)q^{r_1+n}}{\vdots}{\sum_{n\geq0}c_d(n)q^{r_d+n}}\in\mc{H}(\rho)$$
has rational Fourier coefficients. Then there are unique $M_j\in\m$ such that $F=\sum_{j=1}^dM_jD^{j-1}F_0$, and we need to show that in
fact $M_j\in\mq$ for each $j$. A simple inductive argument shows that
$$D^jf_i=\sum_{n\geq0}\beta_{ij}(n)q^{r_i+n},$$
where
$$\beta_{ij}(0)=\prod_{k=0}^{j-1}\left(r_i-\frac{k_0+2k}{12}\right),$$
and $k_0\in\bb{Z}$ denotes the weight of $F_0$. Writing $M_j=\sum_{n\geq0}\alpha_j(n)q^n$ for each $j$, we obtain the formula
$$\cvec{c_1(0)}{\vdots}{c_d(0)}=B\cvec{\alpha_1(0)}{\vdots}{\alpha_d(0)}$$
for the leading Fourier coefficients of $F$; here $B$ denotes the $d\times d$ matrix whose $(i,j)$ entry is $\beta_{ij}(0)$.
Noting that each of these entries is rational, as are the $c_j(0)$, one observes that the invertibility of $B$ would imply $\alpha_j(0)\in\Q$ for each $j$. Now $\beta_{i1}(0)=1$ for each $i$, and for $1<k\leq d-1$ there are polynomials $p_k\in\Q[k_0]$
such that
$$\beta_{ij}(0)=r_i^j+\sum_{k=1}^{j-1}p_kr^{j-k},$$
thus $B$ reduces to the $d\times d$ Vandermonde matrix $V_d$, whose
$(i,j)$ entry is $r_i^{j-1}$. It is well-known that
$$\det(V_d)=\prod_{1\leq i<j\leq d}(r_j-r_i),$$
and by \cite[thm 3]{MM} the $r_j$ are distinct, so $\det(V_d)=\det(B)\neq0$ and we have $\alpha_j(0)\in\Q$ for each $j$. Continuing in this way, one arrives at the formula
$$\cvec{c_1(n)}{\vdots}{c_d(n)}=B\cvec{\alpha_1(n)}{\vdots}{\alpha_d(n)}+\vec{v}_n,$$
for the $n^{th}$ Fourier coefficients of $F$, where the $i^{th}$ entry of $\vec{v}_n$ is a $\Q$-linear combination of the $\alpha_j(k)$,
$1\leq j\leq d-1$, $0\leq k\leq n-1$. Assuming inductively that these entries are rational shows that the $\alpha_j(n)$
are also, and the lemma is proved.\qed

For a prime number $p$, we write $\mc{B}_p\leq\hrhoq$ for the $\mq$-submodule of \emph{$p$-bounded} vectors in $\hrhoq$, \ie the vectors which are not \pubd. Note that if $F_0\in\mc{B}_p$, then clearly $D^jF_0\in\mc{B}_p$ for any $j\geq0$, so by the previous lemma we have
\begin{cor} If $F_0\in\mc{B}_p$, then $\mc{B}_p=\hrhoq$ is a free $\mq$-module of rank $d$.\qed\end{cor}
On the other hand, we have
\begin{prop}\label{prop:ubdmodule}
Suppose there is a prime $p$ such that the Fourier coefficients of the first component of \eqr{F0} satisfy \eqr{ord}.  Then $\mc{B}_p=\{0\}$.
\end{prop}

\proof Let $0\neq g=\sum_{n\geq0}\alpha(n)q^n\in\mq$, so that
$$gf_1=\sum_{n\geq0}\beta(n)q^{r_1+n},\ \ \beta(n)=\sum_{j=0}^n\alpha(j)a(n-j).$$
For any given $n\geq0$, we have by \eqr{ord} that
\begin{equation}\label{eq:beta}
\nu_p(\beta(n))=\nu_p(\alpha(0))+n\delta-\nu_p\left(\prod_{k=1}^nk\lambda(k)\right)
\end{equation}
so long as $\nu_p(\alpha(0)a(n))<\nu_p(\alpha(j)a(n-j))$ for all $1\leq j\leq n$. Again using \eqr{ord}, it is seen that this inequality will hold if
$$(n-j)\delta-\nu_p\left(\prod_{k=n-j+1}^nk\lambda(k)\right)<\nu_p(\alpha(j))-\nu_p(\alpha(0))$$
for each $1\leq j\leq n$. Since $g\in\mq$, there is an integer $M$ such that $\nu_p(\alpha(k))\geq M$ for all $k\geq0$, so for any integer $m\geq0$ satisfying  $m>\nu_p(\alpha(0))-M$, by setting $n=p^m$ in \eqr{beta} we obtain
$$(n-j)\delta-\nu_p\left(\prod_{k=n-j+1}^nk\lambda(k)\right)\leq -m<M-\nu_p(\alpha(0))\leq \nu_p(\alpha(j))-\nu_p(\alpha(0)),$$
for any $1\leq j\leq p^m$. Thus for any such $m$ we have $\nu_p(\beta(n))<-m$, so  $\lim_{m\rightarrow\infty}\nu_p(\beta(p^m))=-\infty$. Since $g$ was arbitrary, we have $\mq F_0\cap\mc{B}_p=\{0\}$.

Now suppose $F\in\mc{B}_p$. Since $D^jF\in\mc{B}_p$ for any $j\geq0$, we obtain from lemma \ref{lem:qform} a relation
$$\cvec{F}{\vdots}{D^{d-1}F}=A\cvec{F_0}{\vdots}{D^{d-1}F_0},$$
with $A=(\alpha_{ij})\in Mat_d(\mq)$. Now, if $A$ were invertible, then we could write $A^{-1}=\det(A)^{-1}C$, where
$\det(A)$ and the entries of the cofactor matrix $C$ lie in $\mq$. But this would yield a relation
$$C\cvec{F}{\vdots}{D^{d-1}F}=\det(A)\cvec{F_0}{\vdots}{D^{d-1}F_0},$$
whose left hand side lies in $\mc{B}_p$, and whose right hand side, by the work of the previous paragraph, does not. This contradiction implies that $A$ is therefore \emph{not} invertible, and consequently there is a relation
$$M_1D^{d-1}F+\cdots+M_dF=0,$$
where at least one of the $M_j\in\mq$ is nonzero. Thus each of the $d$ components of $F$ satisfies the same Fuchsian differential equation of order less than $d$. In particular, the $d$ components must be linearly dependent, so the irreducibility of $\rho$ forces $F=0$, and the proof is complete.\qed

We are now able to complete the proof of theorem \ref{thm:main}. Consider once again an irreducible, three-dimensional representation \eqr{rho3d} with $\rho(T)$ as in \eqr{rhoT3}. By \cite[thm 2.9]{TW}, the eigenvalues $\e{\frac{A}{N}},\e{\frac{B}{N}},\e{\frac{C}{N}}$ of $\rho(T)$ uniquely determine the equivalence class of irreducible representations to which $\rho$ belongs. If the level $N$ of $\rho$ satisfies the hypothesis of corollary \ref{cor:ubd}, then the minimal weight vector \eqr{F0} of $\mc{H}(\rho)$ is \pubd\ for some prime $p$ dividing $N$, thus by proposition \ref{prop:ubdmodule} every nonzero $F\in\hrhoq$ is \pubd. On the other hand, it is clear that the hypothesis of corollary \ref{cor:ubd} will be satisfied by all but a finite number of positive integers $N$, and for each such $N$ there are only a finite number of triples $(A,B,C)$ satisfying the conditions \eqr{ABC}. Thus there are only a finite number of equivalence classes of finite level, irreducible representations $\rho:\G\rightarrow\gln{3}$ that admit vector-valued modular forms with bounded denominators, and theorem \ref{thm:main} is proved.\qed

\section{Finite image representations}\label{sec:monomial}
In this final section, we consider theorem \ref{thm:main} in the classical setting. Here the representation $\rho$ has finite image and, accordingly, the components of its associated vector-valued modular forms are modular for the finite index subgroup $\ker\rho\leq\G$. As we now demonstrate, the application of theorem \ref{thm:main} to this setting yields infinite families of modular forms with unbounded denominators, in support of conjecture \ref{conj:ubd}.

We will utilize \cite[thm 2.1]{RT}, which gives a complete determination of finite image, irreducible representations $\rho:\G\rightarrow\gln{3}$. The classification of \cite{RT} groups the representation classes into \emph{primitive} and \emph{imprimitive} types; in the present setting, $\rho$ is imprimitive if and only if it is monomial, and otherwise is primitive. We first discuss the primitive setting. 

As mentioned in section \ref{sec:background}, every representation of level less than six has a congruence subgroup as kernel, and according to \lc\ the only additional primitive cases yielding finite image are of level seven. From that article and the discussion in section \ref{sec:3d} it follows that four equivalence classes appear here, with $\{A,B,C\}$ in \eqr{ABC} equal to one of $\{0,1,6\},\{0,3,4\},\{1,2,4\},\{3,5,6\}$. In any event, corollary \ref{cor:ubd} does not apply in this situation so we will have nothing further to say about these representations, other than to note that \cite{RT} implies that the first two in the above list are infinite image (thus noncongruence) whereas one may verify that the kernel of each of the last two is congruence of level 7.

Studying the imprimitive representations with finite image amounts to assuming that $\rho=Ind_G^\G(\chi)$ is induced from a finite image character (\ie one-dimensional representation) $\chi:G\rightarrow\cstar$ of an index three subgroup $G\leq\G$. There are only four subgroups of index three in $\G$: $\G^3$, which is the normal subgroup generated by $\{\gamma^3\mid\gamma\in\G\}$, and the conjugate subgroups $(ST)^j\G_0(2)(ST)^{-j}$, $j=0,1,2$, where $\G_0(2)$ denotes the subgroup of matrices in $\G$ that are upper-triangular$\pmod{2}$.

We consider first the subgroup $\G^3$. It is known \cite[pg 36]{Rankin} that the commutator subgroup $(\G^3)'$ is congruence of level 12, and is normal in $\G$. Thus for any character $\chi:\G^3\rightarrow\cstar$ we have that $\ker\chi\geq(\G^3)'$ (since $\cstar$ is an abelian group), and furthermore that the kernel of the representation of $\G$ induced from $\chi$ is equal to $\ker\chi$ (since the intersection of the conjugates of $\ker\chi$ in $\G$ is equal to $\ker\chi$). Thus the vector-valued modular forms associated to such representations have bounded denominators. 

Finally, we consider the three conjugates of $\G_0(2)$. Since inducing characters from conjugate subgroups yields equivalent representations, it is sufficient to consider only the characters of $G=\G_0(2)$. This group has the two cusps $\infty$ and 0, with stabilizers generated by $\pm T$, $\pm ST^2S^{-1}$ respectively, and a single elliptic point $\frac{i-1}{2}$, with stabilizer generated by $(ST)S(ST)^{-1}$. This yields the Fuchsian presentation
$$G\cong\langle E,P_1,P_2\mid E^4=EP_1P_2=1\rangle,$$
with the identifications
\begin{eqnarray*}
E&\leftrightarrow&(ST)S(ST)^{-1},\\
P_1&\leftrightarrow&ST^2S^{-1},\\
P_2&\leftrightarrow&T
\end{eqnarray*}
giving an isomorphism. Thus a character $\chi$ of $G$ is determined \eg\ by any choice of integer $0\leq x\leq3$ and $c\in\cstar$, so that $\chi(E)=\e{\frac{x}{4}}$, $\chi(P_2)=c$. Using the coset decomposition $\G=\cup_{j=0}^2G(ST)^j$, we obtain the induced representation $\rho=Ind_G^\G(\chi)$, with $\rho(\gamma)_{i,j}=\chi((ST)^{i-1}\gamma(ST)^{1-j})$ for each $\gamma\in\G$ (as is customary, here we extend the definition of $\chi$ so that $\chi(g)=0$ for any $g\notin G$). In particular we have
$$\rho(T)=\threemat{\chi(T)}{0}{0}{0}{0}{\chi((ST)T(ST)^{-2})}{0}{\chi((ST)^2T(ST)^{-1})}{0},$$
and from this it follows that the eigenvalues of $\rho(T)$ are $\chi(P_2)=c$ and the two square roots of $\chi(P_1)=c^{-1}\e{-\frac{x}{4}}$.

Assume $\chi$ has finite image, so that the first eigenvalue is $\lambda_1=c=\e{\frac{A}{M}}$, with $0\leq A<M$ and $(A,M)=1$. Then the other eigenvalues are
$$\lambda_2=\e{-\frac{4A+Mx}{8M}}=-\lambda_3,$$
and the level of $\rho$ is $N:=\frac{8M}{(4,Mx)}$ which, since $M\geq1$ is arbitrary, may be any positive even integer. Since $\chi$ has finite image, so does $\rho$, and we have $\ker\rho\leq\ker\chi\leq G$. If $\rho$ is irreducible (generically this will be the case), then (again by \cite[thm 2.9]{TW}) we may assume that $\rho$ is of the form analyzed in section \ref{sec:3d}, and from this and corollary \ref{cor:ubd} we obtain
\begin{prop}
Suppose $\rho:\G\rightarrow\gln{3}$ is as in \eqr{rho3d}, such that the integers $A,B,C,N$ in \eqr{ABC} satisfy the additional constraints $N=2M$, $C=B+M$ for some $M\geq2$. Then $\rho$ has finite image. If, furthermore, there is a prime $p$ satisfying the hypothesis of corollary \ref{cor:ubd}, then the corresponding component of \eqr{F} is a modular form for the noncongruence subgroup $\ker\rho$, and is $\pubd$.\qed
\end{prop}
This result gives infinitely many new examples of noncongruence modular forms with unbounded denominators.

\subsection*{Acknowledgement} The research and preparation involved with this article was conducted while the author was a postdoctoral fellow at the University of Alberta, with funding provided by the Natural Sciences and Engineering Research Council of Canada (NSERC) and the Pacific Institute for the Mathematical Sciences (PIMS). Many thanks to the University -- and in particular to Terry Gannon and Charles Doran -- for support and encouragement.

\bibliographystyle{amsplain}
\bibliography{UBDbib}
\end{document}